\documentclass[12pt,a4paper]{amsart}

\usepackage{amsfonts, amsthm}
\usepackage{enumerate}
\usepackage[english]{babel}
\usepackage{amsmath}
\usepackage{amssymb}
\usepackage{array}
\usepackage[all]{xy}

\swapnumbers

\theoremstyle{theorem}
\newtheorem{theorem}{Theorem}[section]
\newtheorem{proposition}[theorem]{Proposition}
\newtheorem{lemma}[theorem]{Lemma}

\theoremstyle{definition}
\newtheorem{definition}[theorem]{Definition}

\newcommand{\aaa}{\ensuremath{\mathcal{A}}}
\newcommand{\bbb}{\ensuremath{\mathcal{B}}}

\newcommand{\zz}{\ensuremath{\mathfrak{Z}}}

\newcommand{\app}{\ensuremath{\mathsf{App}}}
\newcommand{\PM}{\ensuremath{\mathsf{ProbMet_{}}}}

\newcommand{\uap}{\ensuremath{\mathsf{UAp}}}

\newcommand{\ug}{\ensuremath{\mathsf{UG}}}

\newcommand{\N}{\ensuremath{\mathbb{N}}}

\renewcommand{\inf}{\myinf}
\renewcommand{\min}{\mymin}

\newcommand{\cat}[1]{{\normalfont\mathsf{#1}}}

\DeclareMathOperator*{\myinf}{in\vphantom{p}f}

\DeclareMathOperator*{\mymin}{mi\vphantom{p}n}

\renewcommand{\a}{\ensuremath{\alpha}}
\renewcommand{\b}{\ensuremath{\beta}}

\renewcommand{\d}{\ensuremath{\delta}}

\newcommand{\ve}{\ensuremath{\varepsilon}}

\renewcommand{\l}{\ensuremath{\lambda}}

\def\rw{\rightarrow}

\usepackage{pdfsync}
\begin{document}

\title{A characterisation of probabilistic metrizability for approach spaces}
\author{E. Colebunders, R. Lowen}
\date{}
\maketitle

\vspace{.5cm}

\begin{center}
Abstract
\end{center}

\vspace{.1cm}

 \noindent \scriptsize{Characterisations of metrizable topological spaces or metrizable uniform spaces are well known.  A natural counterpart to being metrizable for topological spaces can be expressed in terms of probabilistic metrizability for approach spaces. The notion of a probabilistic metrizable approach space is based on a well known concrete functor $\Gamma,$ as introduced in \cite{LSY}, from the category of probabilistic metric spaces with respect to a continuous arbitrary t-norm to the category of approach spaces. A characterization of those probabilistic metrizable approach spaces is still missing and in the first part of this paper we solve this problem. 

A natural counterpart to being metrizable for uniform spaces can be expressed in terms of probabilistic metrizability for uniform gauge spaces. In the second part of the paper we start from another concrete functor $\Lambda,$ as described in \cite{CL}, on the category of probabilistic metric spaces with respect to a continuous t-norm to the category of uniform gauge spaces. In a similar way as for the functor $\Gamma$ we obtain a characterisation of probabilistic metrizability of uniform gauge spaces.

The last section of the paper is devoted to an isomorphic description of the category of probabilistic metric spaces. This problem is not new. Previous attempts in providing isomorphic descriptions of the category of probabilistic metric spaces worked with collections of (pseudo)metrics. 
These attempts were only formulated in restricted cases.
 Our isomorphic description is in terms of objects that are sets endowed with a collection of distances, where
  the distances involved do not satisfy the triangle inequality but fulfil a mixed triangle condition instead.

\vspace{.8cm}

 

\section{Introduction}
For many applications the context of approach spaces with contractions is quite suitable as was shown in the context of probability theory \cite{BLC11} and \cite{BHSC}, hyperspaces \cite{LS}, functional analysis \cite{SV} or complexity analysis \cite{CWL} and \cite{CWS}.
A first lax-algebraic description of approach spaces was established by Clementino and Hofmann in \cite{CH}. 
The isomorphic description of $\app$ that was obtained there, became an important example in the development of monoidal topology \cite{MT}, \cite{MMECWT}. In \cite{LSY} and \cite{CL} it became clear that linking a probabilistic metric space to an associated approach space is also quite fruitful. In this paper we pursue this path.

Recently the authors of  \cite{LSY} described a concrete functor $\Gamma: \PM \rw \app,$ from probabilistic metric spaces with non-expansive maps, to the category of approach spaces with contractions \cite{LOW}, \cite{apox}. From \cite{LSY} it is clear that being probabilistic metrizable is a natural counterpart to being metrizable in topological spaces. As in \cite{LSY} our setting will be probabilistic metric spaces with respect to some continuous t-norm $*$ on $[0, 1].$ The image through the functor $\Gamma$ consists of those approach spaces that are probabilistic metrizable. Using an expression from \cite{LSY} the objects in the image of $\PM$ by $\Gamma$ will be described as sets endowed with a local system.
  In the first part of this paper we characterise the image of $\Gamma$ in terms of 
local distances that do not satisfy a triangle inequality but fulfil a mixed triangle condition instead.

Every metric space has an underlying uniform structure and the associated theory of metrizable uniform spaces is well known. As a counterpart we will show that there is a natural concrete functor $\Lambda: \PM \rw \ug$ from probabilistic metric spaces with respect to some continuous t-norm $*$ on $[0, 1]$ with non-expansive maps to the category $\ug$ with objects defined as sets endowed with uniform systems and with uniform contractions as morphisms. The objects in the image of $\PM$ through $\Lambda$ consists of those spaces that are probabilistic metrizable and will be described as sets endowed with a basis for the uniform system consisting of distances. As in the local case, the distances involved do not satisfy the triangle inequality but fulfil a mixed triangle condition instead.

The last section of the paper is devoted to an isomorphic description of the category of probabilistic metric spaces. This problem is not new. Previous attempts in providing isomorphic descriptions of the category of probabilistic metric spaces worked with collections of (pseudo)metrics. These attempts were only formulated in restricted cases, for the minimum t-norm M, or for a weaker form of probabilistic metric spaces, see for instance \cite{N} or \cite{SS}. Our isomorphic desription is in terms of objects that are sets endowed with a collection of distances labeled by the unit interval $]0,1]$ and by morphisms that are levelwise non-expansive. Again the distances involved do not satisfy the triangle inequality but fulfil a mixed triangle condition instead.

\section{Preliminaries}

For more information on probabilistic metric spaces we refer to \cite{SS} and for the link with approach spaces to \cite{LSY} and \cite{CL}.

A function $d: X \times X \rw [0, \infty]$  that is zero on the diagonal  is called a \emph{distance.} This notion is not to be confused with distances between points and sets which also play an important role in approach theory but will not be encountered here. A distance that is symmetric is called a \emph{symmetric distance.}

A \emph{generalised metric} (quasi-metric in \cite{LOW}) on a set $X$ is a map $d: X \times X \rw [0, \infty]$ which is zero on the diagonal and satisfies the triangle inequality. If moreover $d$ is symmetric we call it a \emph{symmetric generalised metric}.

We recall some terminology from \cite{SS} and \cite{LSY}. 

A function $\varphi: [0,\infty] \rw [0,1]$ is called a \emph{distance distribution} if $\varphi$ is monotone,  $\varphi(0) = 0, \varphi(\infty) =1$ and $\varphi$ is left-continuous on $]0,\infty[$.

A binary operation $*$ on the interval $[0,1]$ is a \emph{continuous t-norm} if $([0,1], *, 1)$ is a commutative monoid, which is continuous as a function on $[0,1]^2$ to $[0,1]$ with respect to the usual topologies and satisfies $p * q \leq p' * q'$ whenever $p \leq p'$ and $q \leq q'$ in $[0,1].$

A \emph{probabilistic metric space} is a set $X$ endowed with a map
\begin{equation}\label{PrMet}
\a: X \times X \times [0,\infty] \rw [0,1]
\end{equation}
and a continuous t-norm $*$ such that
\begin{enumerate}
\item [{(P1)}] $\a(x,y,-): [0,\infty] \rw [0,1] $ is a distance distribution,
\item [{(P2)}] $\a (x,x,- ) = \ve_0,$ the largest distance distribution,
\item [{(P3)}] $\a(x,y,r) = \a(y,x,r),$
\item [{(P4)}] $\a(x,y,-) = \ve_0 \Rightarrow  x = y,$
\item [{(P5)}] $\a(y,z,r) * \a(x,y,s) \leq \a(x,z,r+s),$
\end{enumerate}
for all $x,y,z \in X, r,s \in [0, \infty].$

A map $f:(X,\a) \rw (Y,\b)$ between probabilistic metric spaces is \emph{non-expansive} if 
\begin{equation}\label{nonex}
\a(x,x',t) \leq \b(f(x), f(x'),t),
\end{equation}
for all $x,x' \in X, t \in [0, \infty].$

\emph{In this paper we will always assume that the probabilistic metric spaces are defined with respect to some continuous t-norm  which we denote by $*$.}

The category of probabilistic metric spaces with non-expansive maps is denoted by $\PM.$

Next we recall some notions from approach theory. For more details on concepts and results on approach spaces we refer to \cite{LOW} or \cite{apox}.  We recall terminology and basic results that will be needed in this paper.

\begin{definition}\label{localdist}
A collection of ideals $(\aaa(x))_{x\in X}$ in $[0,\infty]^X$ indexed by the points of $X$ is a \emph{local system} if for all $x \in X$ the following properties hold:
\begin{enumerate}
\item[(A1)] $\forall \varphi \in \aaa(x), \ \varphi(x) = 0.$
\item [(A2)] If $\forall \ve > 0, \forall \omega < \infty, \exists \varphi_{\ve}^{\omega} \in \aaa(x)$ such that $ \ \varphi \wedge \omega \leq \varphi_{\ve}^{\omega} + \ve$,\\
 then $\varphi \in \aaa(x).$

\item [(A3)]$\forall \varphi \in \aaa(x), \forall \ve > 0,  \forall \omega < \infty, \exists (\varphi_z)_{z\in x} \in \prod_{z \in X} \aaa(z) \ \text{such that}$
$$
\forall z,y \in X, \varphi (z) \wedge \omega \leq \varphi_x(y) + \varphi_y(z) + \ve.
$$
\end{enumerate}
\end{definition}

The members of $\aaa(x)$ are called \emph{local distances.}

\begin{definition}\label{baselocaldist}
A collection of ideal bases $(\bbb(x))_{x\in X}$ in $[0,\infty]^X$ indexed by the points of $X$ is a \emph{local basis} if for all $x \in X$ the properties  (A1) and (A3) hold.
\end{definition}

In order to derive the local system from a local basis we will require the following \emph{saturation operation}.
Given a subset $\bbb \subseteq [0,\infty]^X$ we define 
\begin{equation}
\langle \bbb \rangle = \{ \varphi \in [0,\infty]^X |\  \forall \ve > 0, \forall \omega < \infty, \exists \varphi_{\ve}^{\omega} \in \bbb \ \text{such that} \ \varphi \wedge \omega \leq \varphi_{\ve}^{\omega} + \ve \},
\end{equation}
the saturation of $\bbb$ and we say that $\langle \bbb \rangle$ is \emph{generated} by $\bbb$.

A map $f: (X, (\aaa(x))_{x\in X}) \rw (Y, (\aaa'(y))_{y\in Y}))$ is a \emph{contraction} if for every $x \in X$ and every $\varphi' \in \aaa'(f(x))$ we have
\begin{equation}\label{contractsys}
\varphi' \circ f \in \aaa(x).
\end{equation}
It is sufficient to impose the condition on a local basis in the codomain.
The category $\cat{App}$ of approach spaces has as objects sets endowed with a local system and contractions as morphisms.
For other \emph{isomorphic descriptions} of $\cat{App}$ we refer to the literature.

\section{Probabilistic metrizable approach spaces}

In this section we give an internal characterisation of probabilistic metrizable approach spaces.
We start with the concrete functor  \cite{LSY}, \cite{CL}
\begin{equation}\label{faith}
\Gamma: \PM \rw \app,
\end{equation}
where to each probabilistic metric space $(X, \a, *)$ an approach space is assigned.
In \cite{LSY} the functor was described using approach spaces in terms of distances $\d(x,A)$  between points and sets. In \cite{CL} the functor was equivalently described in terms of the tower of an approach space and it was shown that the codomain of the functor can be restricted to $\uap$, the subcategory of $\app$ of uniform approach spaces.
In \cite{LSY} yet another equivalent transition from probabilistic metric spaces $(X, \a, *)$ to approach spaces in terms of local systems was formulated using countable collections of local distances.
It is a similar formulation of the functor $\Gamma$ we will use here.

Let $(X,\a,*)$ be a probabilistic metric space. In \ref{uap2} we show that an \emph{underlying approach space} $X$ is associated with it, defined by the following local basis in $x,$
\begin{equation}\label{associated}
\bbb(x) = \{ \varphi_{\l,x} \ | \ 0 < \l \leq 1 \} \ \ 
\end{equation}
with 
$$  \varphi_{\l,x} (y) = \inf \{ 0 \leq \gamma < \infty \ | \ \a(x,y,\gamma) > 1 - \l \}.$$

\begin{proposition}\label{}
For $0 \leq \gamma < \infty$ and for $0 < \l \leq 1$ we have
\begin{equation}\label{equivlocal}
\varphi_{\l,x} (y) < \gamma \Leftrightarrow \a(x,y,\gamma) > 1-\l.
\end{equation}
\begin{proof}
Suppose $\varphi_{\l,x} (y) < \gamma,$ then there exists $\mu < \gamma$ satisfying $\a(x,y,\mu) > 1-\l.$
This implies $\a(x,y,\gamma) > 1-\l$ since $\a$ is nondecreasing.

To prove the other implication, assume $\a(x,y,\gamma) > 1-\l.$ Left continuity of $\a$ implies that there exists $\mu < \gamma$ with $ \a(x,y,\mu) > 1-\l.$ This implies\\
 $\varphi_{\l,x} (y) \leq \mu < \gamma.$
\end{proof}
\end{proposition}

\begin{proposition}\label{uap2}
$\bbb(x) = \{ \varphi_{\l,x} \ | \ 0 < \l \leq 1 \} $ as defined in \eqref{associated} is a local basis on $X$ in the sense of \ref{baselocaldist}.
\begin{proof}
Clearly all $\varphi_{\l,x} \in [0,\infty]^X.$
That (A1) is fulfilled, namely $\varphi_{\l,x} (x) = 0$ for every $\l,x$ follows from the fact that $\a(x,x,\gamma) = 1$ for all $\gamma>0$.

Next we prove condition (A3). In fact we prove a \emph{a local mixed triangle inequality} {\bf(LM)} 
that clearly implies (A3):
 For $0 < \ve \leq 1$ and
$0<\l \leq 1$ and such that 
$$
(1-\l) * (1-\l) > 1 - \ve.
$$
and for all $x,y,z \in X$ we have
\begin{equation}\label{exists}
\varphi_{\ve,x} (z) \leq \varphi_{\l,x} (y) + \varphi_{\l, y} (z).
\end{equation}
First observe that in view of the continuity of the t-norm $*$ on $[0,1]$ for every $0 < \ve \leq 1$ there exists $0<\l \leq 1$ such that 
$$
(1-\l) * (1-\l) > 1 - \ve.
$$

If the righthandside of \eqref{exists} is $\infty$ there is nothing to prove. 
Next let $\gamma, \gamma'$ be arbitrary, with $\varphi_{\l,x} (y) < \gamma$ and $\varphi_{\l, y} (z) <\gamma'.$
Applying \eqref{equivlocal} we get $\a(x,y,\gamma) > 1 - \l$ and $\a(y,z,\gamma') > 1 - \l$.
By axiom (P5) for $\a$ we obtain 
$$
\a(x,z, \gamma + \gamma') \geq \a(y,z,\gamma') * \a(x,y,\gamma) \geq (1-\l) * (1-\l) > 1 - \ve.
$$
Again applying \eqref{equivlocal} we get $\varphi _{\ve,x}(z) < \gamma+ \gamma',$ from which the conclusion \eqref{exists} follows.

\end{proof}
\end{proposition}

It is easy to see that $\Gamma$ is indeed a functor by applying the results in \cite{LSY}. This fact also follows from \cite{CL} as the description given there in terms of the tower is equivalent to the one presented here, since the sets $\{\varphi _{\l,x} < \rho\}$ coincide with the sets $V_{\l,x}^\rho$ \cite{CL}.
Nevertheless we present an explicit proof as we will need the details further on.

\begin{proposition}\label{functoruap}
$\Gamma: \PM \rw \app$ is a concrete functor, meaning that when  $f: (X,\a) \rw (Y, \a')\  \text{is non-expansive}$ then $\Gamma(f) :  \Gamma((X, \a)) \rw \Gamma((Y,\a'))$ is a contraction, where $\Gamma((X, \a)) $ and $\Gamma((Y, \a')) $ have respective local systems
$(\aaa(x))_x$ generated by $\bbb(x) =\{ \varphi_{\l,x} \ | 0 < \l \leq 1\} $ and $(\aaa'(y))_{y}$ generated by $\bbb'(y)=\{ \varphi'_{\l,y} \ | \ 0 < \l \leq 1\} $ as in \eqref{associated}.
\begin{proof}
First observe that in case $\varphi'_{\l,f(x)} (f(x'))= \infty$ then for every $\gamma < \infty,$\\
$ \a(f(x), f(x'), \gamma) \not > 1 - \l.$
This implies that for every $\gamma < \infty,\  \a(x,x',\gamma) \not > 1 - \l.$ Hence we can conclude that $\varphi _{\l, x}(x') = \infty.$

 In the finite case, applying \eqref{equivlocal} we obtain the following equivalences:
 
 \begin{eqnarray*}
  f: (X,\a) \rw (Y, \a')\  \text{is non-expansive} \ &\Leftrightarrow & \forall t, \ \a(x,x',t) \leq \a'(f(x), f(x'), t )\\
&\Leftrightarrow& \forall t, \forall \l : \ \ (\a(x,x', t) > 1 - \l  \Rightarrow\\
&& \a'(f(x), f(x') , t) > 1 - \l )\\
&\Leftrightarrow& \forall t,  \forall \l :  \  \varphi_{\l,x} (x') < t \Rightarrow \varphi ' _{\l, f(x)} (f(x')) < t\\
&\Leftrightarrow& \forall \l, \  \varphi'_{\l,f(x)} (f(x')) \leq \varphi _{\l,x} (x').\\
\end{eqnarray*}
Clearly the last line implies:
$$
\text{if} 
\  \varphi' _{\l, f(x)} \in \bbb'(f(x))\   \text{then}\ \ \varphi'_{\l, f(x)} \circ f \in \aaa(x),
$$
and by \eqref{contractsys} $\Gamma(f)$ is a contraction.
\end{proof}
\end{proposition}

In order to determine the image of $\PM$ by $\Gamma$ we prove the following properties.

\begin{theorem}[\bf Part I of Theorem 3.8]\label{char}
 Let $(X,\a,*)$ be a probabilistic metric space. The associated local basis $\bbb(x) = \{ \varphi_{\l,x} \ | \ 0 < \l \leq 1 \} $ in $x,$ as in \eqref{associated} and 
 in \ref{uap2}, has the following properties:
\begin{enumerate}
\item [\bf{(LS)}] Symmetry: For every $0 < \l \leq 1$ and for $x,y \in X$, we have 
$$\varphi_{\l,x} (y) = \varphi _{\l,y} (x).$$
\item [\bf{(LD)}] Density:  For every $x$ and $0 < \l \leq 1$ 
$$
\varphi_{\l,x} = \inf _{\rho < \l} \varphi_{\rho,x}.
$$ 
\item[\bf{(LT)}] t-Norm: For $0 < \ve \leq 1,$ \  $0 < \l \leq 1$ and $0 < \l' \leq 1$
with
$$
 (1 - \l') * ( 1 - \l) > 1- \ve
$$
then for all $x,y,z \in X$ we have
$$
\varphi_{\ve,x} (z) \leq \varphi _{\l, x} (y) + \varphi _{\l', y} (z).
$$
\item [\bf{(LH)}] Separation:
$$
x \not = y \ \text{then} \  \exists\l, \  0 < \l \leq 1, \ \varphi_{ \l,x} > 0.
$$
\end{enumerate}

\begin{proof}
{\bf{(LS)}} Follows from (P3), since $\a(x,y,\gamma) = \a(y,x,\gamma)$ for all $x,y,\gamma.$\\
{\bf(LD)} Suppose $0 < \l' \leq \l \leq 1$ and $x \in X.$ First observe that in the infinite case, when $\varphi_{\l,x} (y) = \infty$ then $\{ \gamma \ | \ \a(x,y,\gamma) > 1 - \l \} = \emptyset $ and then clearly also $\{ \gamma \ | \ \a(x,y,\gamma) > 1 - \l' \} = \emptyset $ and $\varphi_{\l',x} (y) = \infty.$
In the finite case  it is sufficient to apply \eqref{equivlocal} in order to conclude that $\varphi_{\l,x}  \leq  \varphi_{\l',x}.$\\
From the previous observation we already have $\varphi_{\l,x} \leq \inf _{\rho < \l} \varphi_{\rho,x}.$
We show the reverse inequality. Let $y \in X$. In the infinite case with $\varphi _{\l, x} (y) = \infty,$ there is nothing to show.
Next assume $\varphi_{\l, x} (y) <  \gamma < \infty.$ Applying \eqref{equivlocal} we have $\a(x,y,\gamma) > 1 - \l.$ Then we can choose $\rho < \l$ with 
$\a(x,y,\gamma) > 1 - \rho.$ Clearly $\varphi_{\rho,x} (y) < \gamma$ which implies $ \inf _{\rho < \l} \varphi_{\rho,x}(y) < \gamma.$\\
{\bf(LT)} 
Let $0 < \ve \leq 1,$ \  $0 < \l \leq 1$ and $0 < \l' \leq 1$ be satisfying
$$ (1 - \l') * ( 1 - \l) > 1 - \ve.$$ 
For $x,y,z \in X$ arbitrary,
if the righthandside of the inequality $\varphi_{\ve,x} (z) \leq \varphi _{\l, x} (y) + \varphi _{\l', y} (z)$ is infinite, there is nothing to prove. Next assume $\varphi _{\l, x} (y) < \gamma$ and  $\varphi _{\l', y} (z) < \gamma'$. Then we have $\a(y,z, \gamma') > 1 - \l'$ and $\a(x,y, \gamma) > 1 - \l.$
Now apply (P5)
$$
\a(x,z,\gamma + \gamma') \geq \a(y,z, \gamma')  * \a(x,y, \gamma) \geq (1 - \l') * ( 1 - \l)  > 1-\ve.
$$
Finally we get $\varphi_{\ve,x} (z) < \gamma + \gamma'.$\\
{\bf(LH)} 
Suppose $x \not = y.$ That there exists $0 < \l \leq 1$ with $\varphi_{ \l,x} > 0$ follows from the fact that by (P4) $\exists \gamma > 0$ with $\a(x,y,\gamma) \not = 1$ and hence $\exists \l$ with $\a(x,y,\gamma )\not  > 1 - \l.$
\end{proof}
\end{theorem}

In fact the condition t-Norm {\bf(LT)} clearly implies the local mixed triangle condition {\bf(LM)} of \ref{uap2}, so {\bf(LT)} can be seen as a \emph{local strong mixed triangle inequality.}

Next we investigate the other direction, namely the question how to associate a probabilistic metric space to a particular local basis of an approach space.

\begin{definition}\label{}
Let $X$ be an approach space, having a local basis\\
 in $x$
$$
\bbb(x) = \{ \varphi_{\l,x} \ | \ 0  < \l \leq 1 \}.
$$
For $x,y \in X$ we define $\b (x,y, \infty ) = 1$ and further for $\gamma < \infty$:
\begin{equation}\label{beta}
\b(x,y,\gamma) = \sup \{1 - \l \ | \ \varphi_{\l,x} (y) < \gamma\} = 1 - \inf  \{\l \ | \ \varphi_{\l,x} (y) < \gamma\}.
\end{equation}
\end{definition}

\begin{proposition}\label{equivbeta}
Given an approach space $X$ having a local basis satisfying {\bf(LD)},
then for $0 \leq \gamma < \infty$ we have:
$$
\b(x,y,\gamma) > 1 - \l   \Leftrightarrow  \varphi_{\l,x} (y) < \gamma.
$$
\begin{proof}
Let  $\sup \{1 - \rho \ | \ \varphi_{\rho,x} (y) < \gamma\}  > 1 - \l.$ Then $\exists \rho < \l, \ \varphi_{\rho,x}(y) < \gamma$ and applying {\bf(LD)} we have 
$\varphi_{\l,x}(y) < \gamma.$

For the other implication,  assume $\varphi_{\l,x}(y) < \gamma.$ Then again applying {\bf(LD)}\\
 $\exists \rho < \l, 
 \  \ \varphi_{\rho,x}(y) < \gamma.$
Hence $1 - \l < 1 - \rho \leq \b(x,y,\gamma).$
\end{proof}
\end{proposition}

 First we need the following lemma.
\begin{lemma}\label{star}
Let $([0,1], *)$ be a continuous t-norm. For all $a,b,d \in ]0,1],$ the following are equivalent:
\begin{enumerate}
\item $d \geq a * b$
\item $\forall \l,\l' \in ]0,1], $\  if $a > 1-\l, \ b > 1-\l'$ then $d \geq (1-\l') * (1-\l)$
\item $\forall \rho  \in ]0,1],$ \ if $a * b > 1-\rho$ then $d \geq 1-\rho.$
\end{enumerate}
\begin{proof}
That (1) and (3) are equivalent and that (1) implies (2) are straightforward.
We prove that (2) implies (3). 

Let $ \rho  \in ]0,1],$ and assume that $a*b > 1-\rho.$
Choose strictly increasing sequences $(1-\l_n)_n,$ in $]0,1]$ converging to $a$ and $(1-\l'_n)_n,$ converging to $b.$ 
By the continuity of $*$ on $[0,1]\times [0,1]$ the sequence $((1-\l_n) * (1-\l'_n))_n$ converges to $a * b$ with $(1-\l_n) * (1-\l'_n) \leq a*b$ for all $n.$
Since $]1-\rho, 1]$ is an open neighborhood of $a * b$ we can fix $n \in \N$ such that $(1-\l_n) * (1-\l'_n) \in ]1-\rho, 1].$
By application of (2) we now have
$$
d \geq (1-\l_n) * (1-\l'_n) > 1 - \rho.
$$

\end{proof}
\end{lemma}

\begin{theorem} [\bf Part II of Theorem 3.8]\label{pmetr}
If $X$ is an approach space having a local basis in $x$
$$
\bbb(x) = \{ \varphi_{\l,x} \ | \ 0  < \l \leq 1 \},
$$
satisfying {\bf(LS)}, {\bf(LD)}, {\bf(LT)} for some continuous t-norm $*$ and {\bf(LH)}, then there exists a probabilistic metric space $(X,\b,*)$ with underlying approach space $X$.

\begin{proof}
Suppose 
$X$ has a local basis in $x$ 
$$
\bbb(x) = \{ \varphi_{\l,x} \ | \ 0  < \l \leq 1 \},
$$
satisfying {\bf(LS), (LD), (LT)} and {\bf(LH)}.

On $X$ define $\b(x,y, -)$ as in \eqref{beta}.
We prove that $(X, \b, *) $ is a probabilistic metric space in the sense of \eqref{PrMet}.

$\b(x,y,0) =0$ as $\varphi_{\l,x} (y) \geq 0$ and $\b(x,y,\infty) =1$ by definition.

$\b(x,y,-)$ is nondecreasing as for $\gamma \leq \gamma',$  $\{ \l \ | \ \varphi_{\l,x} (y) < \gamma\} \subseteq \{ \l \ | \ \varphi_{\l,x} (y) < \gamma'\}.$

$\b(x,y,-)$ is left continuous on $]0, \infty[$. That $\sup_{s < \gamma} \b(x,y,s) \leq  \b(x,y,\gamma)$ is clear from the previous property. For the reverse inequality assume $\b(x,y,\gamma) > 1 - \l$. By \ref{equivbeta} this implies that $\varphi _{ \l,x} (y) < \gamma$  and hence $\exists s < \gamma,  \ \ \varphi _{ \l,x} (y) < s.$ This implies that $\exists s <\gamma,  \  \b(x,y,s) > 1 - \l.$ So we have that $\sup_{s < \gamma} \b(x,y,s) > 1 - \l.$
From the previous claims we can conclude that (P1) holds.

To see that (P2) holds,  observe that for $\gamma >0$  we have\\
 $$\b(x,x,\gamma) = \sup \{1 - \l \ | \ \varphi_{\l,x} (x) < \gamma\}  = 1.$$

Clearly (P3) follows immediately from {\bf(LS)}.

Next we check (P4). Suppose $x \not = y,$  from {\bf (LH)} it follows that $$\exists \l, \ \varphi_{\l,x} (y) \not = 0.$$
Choosing $\gamma$ such that $0 < \gamma < \varphi_{\l,x} (y), $ from \ref{equivbeta} we obtain $\b(x,y,\gamma) \not = 1.$

Next we check that $\b$ satisfies (P5). Let $x,y,z \in X,$ and assume $\b(y,z,\gamma') > 1 - \l'$ and $\b(x,y,\gamma) > 1 - \l.$  
For every $ 0 < \ve \leq 1 $ satisfying 
$$
(1 - \l') * ( 1 - \l) > 1-\ve,
$$
applying {\bf(LT)},  
we have
$$
\varphi_{\ve,x} (z) \leq \varphi _{\l, x} (y) + \varphi _{\l', y} (z) < \gamma' + \gamma.
$$
Hence
$$
\b(x,z, \gamma + \gamma') > 1 - \ve. 
$$
This implies
$$
\b(x,z, \gamma + \gamma') \geq  (1 - \l') * ( 1 - \l).
$$
Applying \ref{star}
we can conclude that $\b(x,z, \gamma + \gamma') \geq \b(y,z, \gamma')  * \b(x,y, \gamma). $\\

Finally we prove that the underlying approach space of $(X,\b)$ is $X$.\\
Let $\{\psi_{\l,x} \ | \ 0 < \l \leq 1 \}$ be associated with $(X,\b)$ as in \eqref{associated}. Then from \eqref{equivlocal}  and \ref{equivbeta} for $\gamma$ finite
$$
\psi_{\l,x}(y) < \gamma \Leftrightarrow \b(x,y, \gamma) > 1 - \l \Leftrightarrow \varphi_{\l,x}(y) < \gamma.
$$
Hence we can conclude that \ $\psi_{\l,x} = \varphi_{\l,x}$ for all $\l$. 
\end{proof}
\end{theorem}

From Theorems \ref{char} (Part I of 3.8) and \ref{pmetr}  (Part II of 3.8), we now have the following conclusion.

\begin{theorem}
The probabilistic metrizable approach spaces are exactly those $X$ having a local basis in $x$
$$
\bbb(x) = \{ \varphi_{\l,x} \ | \ 0  < \l \leq 1 \},
$$
satisfying {\bf(LS)}, {\bf(LD)}, {\bf(LT)} for some continuous t-norm $*$ and {\bf(LH)}. 
\end{theorem}

\section{Probabilistic metrizability of uniform gauge spaces}

In this section we give an internal characterisation of probabilistic metrizablility in the uniform case.  We will work with objects defined by a uniform alternative of the local systems, namely uniform systems \cite{LW} and with uniform contractions as morphisms. The category we will be working with is isomorphic to the category $\ug$ of uniform gauge spaces \cite{LOW}. 
We start with the concrete functor  \cite{CL}
\begin{equation}\label{concreteug}
\Lambda: \PM \rw \ug,
\end{equation}
where to each probabilistic metric space $(X, \a, *)$ an underlying space $X$ is assigned endowed with a uniform system.

Uniform systems are defined in terms of collections of distances between points.  That means that the triangle inequality of the functions in $[0,\infty] ^{X \times X}$ is not required. As in the previous section it is replaced by a mixed triangle condition (AU3).
For more details we refer to \cite{LW}.

\begin{definition}\label{usystems}
Let $X$ be a set. A \emph{uniform system} on $X$ is an ideal $\Upsilon \subseteq [0,\infty] ^{X \times X}$ such that 
\begin{enumerate}
\item [(AU1)] $\forall d \in \Upsilon, \forall x \in X, d(x,x) = 0.$
\item [(AU2)] If $\forall \ve >0, \forall N < \infty, \exists d_\ve^N \in \Upsilon, e \wedge N \leq d_\ve^N + \ve, \ \text{then} \ e \in \Upsilon.$
\item [(AU3)] $ \forall d \in \Upsilon, \forall N < \infty, \exists e^N \in \Upsilon, \forall x,y,z \in X,$
$$
d(x,z) \wedge N \leq e ^N (x,y) + e^N (y,z).
$$
\item [(AU4)] $\forall d \in \Upsilon , d^{-1} \in \Upsilon.$
\end{enumerate}
\end{definition}

The members of $\Upsilon$ are called \emph{uniform distances.}

As for the local systems, also in the uniform context it is sufficient to have a basis.
\begin{definition}\label{ubsystem}
Let $X$ be a set and $\Psi \subseteq [0,\infty] ^{X \times X}$. Then $\Psi$ is said to be a \emph{uniform basis} if $\Psi$ is an ideal basis satisfying (AU1), (AU3) and (AU4).
\end{definition}

In order to derive the uniform system from a uniform basis, we will require the following \emph{saturation operation}.
Given a subset $\Psi \subseteq [0,\infty]^{X \times X},$ a basis of a uniform system, we define 
\begin{equation}
\langle \Psi \rangle = \{d \in [0,\infty] ^{X \times X} \ | \ \forall \ve >0, \forall N < \infty, \exists d_\ve^N \in \Psi, e \wedge N \leq d_\ve^N + \ve, \}
\end{equation}
the saturation of $\Psi$ and we say that $\langle \Psi \rangle$ is \emph{generated} by $\Psi$.

A function $f: (X, \Upsilon) \rw (Y, \Upsilon')$ is a \emph{uniform contraction} iff 
\begin{equation}\label{contractsystem}
 d' \circ (f \times f)\in \Upsilon  \ \text{whenever} \ \ d \in \Upsilon'.
\end{equation}
Again it is sufficient to impose the condition on a uniform basis in the codomain.

In \cite{LW} it is shown that the category of sets endowed with uniform systems and uniform contractions is an \emph{isomorphic description} of $\ug,$ the category of uniform gauge spaces as introduced in \cite{LOW}. Hence we keep on calling the category $\ug$.

Next we show that there is a concrete functor
$$
\Lambda : \PM \rw \ug.
$$

In \cite{CL} the object associated to a probabilistic metric space was described in terms of the tower.
In \cite{N} \cite{SS}, yet another equivalent transition from probabilistic metric spaces in terms of collections of (finitely valued) symmetric distances was given.
It is the same formulation for the associated collections of (symmetric) distances we will use here. The difference is that we interpret the collection as a basis for a uniform system. We come back to the setting of collections of (symmetric) distances in the last section.

Let $(X,\a,*)$ be a probabilistic metric space. An \emph{underlying $\ug$-space} $X$ is associated with it and is defined by giving the following uniform basis $\Psi,$
\begin{equation}\label{uassociated}
\Psi = \{ d_\l \ | \ 0 < \l \leq 1 \}, \ \ 
\end{equation}
with 
$$  d_\l(x,y) = \inf \{ 0 \leq \gamma < \infty \ | \ \a(x,y,\gamma) > 1 - \l \}.$$

Quite similar as in \eqref{equivlocal} we can prove that for $\gamma$ finite we have the following equivalence. Because of the similarities with the previous section we do not repeat all the proofs.

\begin{proposition}
For $0 \leq \gamma < \infty$ and for $0 < \l \leq 1$ we have
\begin{equation}\label{equivug}
d_\l (x,y) < \gamma \Leftrightarrow \a(x,y,\gamma) > 1-\l.
\end{equation}
\end{proposition}

\begin{proposition}\label{ug2}
$\Psi = \{ d_\l \ | \ 0 < \l \leq 1 \} $ as defined in \eqref{uassociated} is a uniform basis on $X$ in the sense of \ref{ubsystem}.
\begin{proof}
By definition all $d_\l \in [0,\infty]^{X \times X}.$
That (AU1) is fulfilled, namely $d_\l (x,x) = 0$ for every $\l,$ follows from the fact that $\a(x,x,\gamma) = 1$ for all $\gamma>0$.

Next we prove condition (AU3). In fact we prove a \emph{uniform mixed triangle inequality} {\bf(UM)}
that clearly implies (AU3):
For $0 < \ve \leq 1,$ \ 
$0<\l \leq 1$ such that 
$$
(1-\l) * (1-\l) > 1 - \ve.
$$
and for all $x,y,z \in X$ we have
\begin{equation}\label{ustar}
d_\ve (x,z) \leq d_{\l }(x,y) + d_{\l }(y,z).
\end{equation}
First observe that in view of the continuity of the t-norm $*$ on $[0,1]$ there exists $0<\l \leq 1$ such that 
$$
(1-\l) * (1-\l) > 1 - \ve.
$$
If the righthandside of \eqref{ustar} is $\infty$ there is nothing to prove. 
Next let $\gamma, \gamma'$ be arbitrary, with $d_{\l}(x,y) < \gamma$ and $d_{\l} (y, z) <\gamma'.$
Applying \eqref{equivug} we get $\a(x,y,\gamma) > 1 - \l$ and $\a(y,z,\gamma') > 1 - \l$.
By axiom (P5) for $\a$ we obtain 
$$
\a(x,z, \gamma + \gamma') \geq \a(y,z,\gamma') * \a(x,y,\gamma) \geq (1-\l) * (1-\l) > 1 - \ve.
$$
Again applying \eqref{equivug} we get $ d_\ve (x,z) < \gamma + \gamma',$ from which the conclusion \eqref{ustar} follows.

The property (AU4) follows from condition {\bf(US)} which we list below in \ref{uchar}.

\end{proof}
\end{proposition}

It is easy to see that $\Lambda$ is indeed a functor by applying the results in \cite{CL}, as the description given there in terms of the tower is equivalent to the one presented here. The sets $\{d_\l < \rho\}$ coincide with the sets $U_\l^\rho$ \cite{CL}.
The proof  in terms of uniform distances is quite similar to the one in \ref{functoruap}.

\begin{proposition}\label{uimage}
$\Lambda: \PM \rw \ug$ is a concrete functor, meaning that when $f: (X,\a) \rw (Y, \a')$ is non-expansive then  $\Lambda(f) : (X, \Upsilon) \rw (Y, \Upsilon')$ is a uniform contraction in the sence of \eqref{contractsystem}, with $\Lambda((X,\a)) = (X, \Upsilon)$ and $\Lambda((Y,\a')) = (Y, \Upsilon')$
\end{proposition}

In order to determine the image of $\PM$ by $\Lambda$ we list the following properties.

\begin{theorem}[\bf Part I of Theorem 4.10]\label{uchar}
 Let $(X,\a,*)$ be a probabilistic metric space. The associated uniform basis $\Psi = \{ d_\l \ | \ 0 < \l \leq 1 \}$
 as in \eqref{uassociated} has the following properties: 
\begin{enumerate}
\item [{\bf(US)}] Symmetry: For every $0 < \l \leq 1$ and for $x,y \in X$, we have 
$$d_\l (x, y) = d_\l (y, x).$$
\item [{\bf(UD)}] Density:  For every $0 < \l \leq 1$ 
$$
d_\l = \inf _{\rho < \l} d_\rho.
$$ 
\item[{\bf(UT)}] t-Norm: For $0< \ve \leq 1,$ \ $0 < \l \leq 1$ and $0 < \l' \leq 1$ with 
$$
 (1 - \l') * ( 1 - \l) > 1 - \ve,
$$
 and for $x,y,z \in X$ we have
$$
d_\ve(x, z) \leq d_\l(x, y) + d_{\l'}(y, z).
$$
\item [{\bf(UH)}] Separation:
$$
x \not = y \ \text{then} \  \exists\l, \  0 < \l \leq 1, \ d_\l (x,y) > 0.
$$

\end{enumerate}

\end{theorem}
Similar to the local situation it is easily observed that {\bf(UT)} implies {\bf(UM)}. So {\bf(UT)} can be seen as a \emph{uniform strong mixed triangle condition.}

Next we investigate the other direction, namely the question how to associate a probabilistic metric space to a particular uniform basis of a $\ug$-space.

\begin{definition}\label{reconstruction}
Let $X$ be a $\ug$-space having a uniform basis
$$
\Psi = \{ d_\l \ | \ 0  < \l \leq 1 \}.
$$
On $X$ we define $\b (x,y, \infty ) = 1$ and further for $\gamma < \infty$:
\begin{equation}\label{ubeta}
\b(x,y,\gamma) = \sup \{1 - \l \ | \ d_\l(x, y) < \gamma\} = 1 - \inf  \{\l \ | \ d_\l (x, y) < \gamma\}.
\end{equation}

\end{definition}

\begin{proposition}\label{uequivbeta}
Given $X$ having a uniform basis satisfying {\bf(UD)},
for $0 \leq \gamma < \infty$ we have:
$$
\b(x,y,\gamma) > 1 - \l   \Leftrightarrow  d_\l(x, y) < \gamma.
$$
\end{proposition}

\begin{theorem}[\bf Part II of Theorem 4.10]\label{upmetr}
If $X$ is a $\ug$-space having a uniform basis
$$
\Psi = \{ d_\l \ | \ 0  < \l \leq 1 \},
$$
satisfying {\bf(US), (UD), (UT)} for some continuous t-norm $*$ and {\bf(UH)}, then there exists a probabilistic metric space $(X,\b, *)$ with underlying $\ug$-space $X$.

\begin{proof}
Suppose
$X$ has a uniform basis 
$$
\Psi = \{ d_\l  \ | \ 0  < \l \leq 1 \}
$$
satisfying {\bf(US), (UD), (UT)} and {\bf(UH)}.

On $X$ define $\b(x,y, -)$ as in \eqref{ubeta}.
Similar to the proof in \ref{pmetr},
$(X, \b, *) $ is a probabilistic metric space.

Finally we prove that the underlying $\ug$-space of $(X,\b)$ is $X$.

Let $\{ e_\l \ | \ 0 < \l \leq 1 \}$ be associated with $(X,\b)$ as in \eqref{uassociated}. Then from \eqref{equivug}  and \ref{uequivbeta} for $\gamma<\infty$ 
$$
e_\l(x, y) < \gamma \Leftrightarrow \b(x,y, \gamma) > 1 - \l \Leftrightarrow d_\l(x, y) < \gamma.
$$
Hence we can conclude that $e_\l = d_{\l}$ for all $\l$.
\end{proof}
\end{theorem}

From Theorems \ref{uchar} (Part I of 4.10) and \ref{upmetr} (Part II of 4.10), we now have the following conclusion.

\begin{theorem}
The probabilistic metrizable $\ug$-spaces are exactly those $X$ having a uniform basis 
$$
\Psi = \{ d_\l \ | \ 0  < \l \leq 1 \},
$$
satisfying {\bf(US)}, {\bf(UD)}, {\bf(UT)} for some continuous t-norm $*$ and {\bf(UH)}. 
\end{theorem}

\section{An isomorphic description of the category of Probabilistic Metric Spaces}

Several attempts have been made to describe the category of probabilistic metric spaces with non-expansive maps isomorphically in terms of collections of symmetric generalised metrics with levelwise non-expansive maps.  However in order to obtain an isomorphism based on generalised metrics, one had to make restrictions. In \cite{SS} the continuous t-norm  is restricted to M, taking  $a * b = \min (a,b)$. Another direction of research is to express weak probabilistic metric spaces in terms of symmetric generalised metrics, which is realised by weakening the axiom (P5) and replacing it by the axiom \cite{N} 
$$
\a(x,y,\gamma) =1\ \text{and} \ \a(y,z, \gamma') =1 \Rightarrow \a(x,z, \gamma + \gamma') = 1.
$$

or by imposing 

$$
\a(x,y,\gamma) > t\ \text{ and } \ \a(y,z, \gamma' )> t  \Rightarrow \a(x,z, \gamma + \gamma') >t
$$
for $0 \leq t < 1.$

In this section, as in the previous ones, we stay with the axiom (P5) formulated for an arbitrary continuous t-norm as in \eqref{PrMet} and we work with collections of distances where the triangle condition of the members of the collection is replaced by a mixed triangle condition of the 
collection itself. 

Note that in what follows there is a fundamental difference compared to the previous section. In Definition \ref{Z} the objects will be sets endowed with a collection of distances. Different collections will stand for different objects. However in section 4 a collection of distances is a basis for a uniform system. This implies that different collections of distances may generate the same $\ug$-object. A second fundamental difference lies in the formulation of the morphisms in \eqref{morphisms}, where they are levelwise non-expansive. However in section 4 morphisms are defined as uniform contractions \eqref{contractsystem}, where in the domain of the functions the whole uniform system is involved.

\begin{definition}\label{Z}
Let $\zz$ be the category with objects, $(X, \{d_\l \ | \ 0< \l \leq 1\}),$ sets endowed with a collection of distances numbered by $ 0< \l \leq 1,$ where\\ $\{d_\l \ | \ 0< \l \leq 1\}$
satisfies the conditions: {\bf(US)}, {\bf(UD)}, {\bf(UT)}  for some continuous t-norm $*$  and {\bf(UH)} from \ref{uchar}.

A morphism  $f : (X, \{d_\l \ | \ 0< \l \leq 1\}) \rw (Y, \{d'_\l \ | \ 0< \l \leq 1\})$ is \emph{levelwise non-expansive} if
\begin{equation}\label{morphisms}
\forall \l, \  d'_\l(f(x), f(x')) \leq d_\l (x,x'), \ \text{whenever} \ x, x' \in X.
\end{equation}
\end{definition}

\vspace{.7cm}

Let $(X,\a,*)$ be a probabilistic metric space. As in \cite{N}
\begin{equation}\label{distanceassociated}
(X, \{ d_\l \ | \ 0 < \l \leq 1 \})  
\end{equation}
 is associated with $(X,\a,*)$ by
$$  d_\l(x,y) = \inf \{ 0 \leq \gamma < \infty \ | \ \a(x,y,\gamma) > 1 - \l \}.$$
This is the same formula as in \eqref{uassociated} and again we have the equivalences of \eqref{equivug}, namely 
for $0 \leq \gamma < \infty$ and for $0 < \l \leq 1$ 
 \begin{equation}\label{isoequiv}
d_\l (x,y) < \gamma \Leftrightarrow \a(x,y,\gamma) > 1-\l.
\end{equation}

\begin{proposition}\label{delta}
This defines a concrete functor $\Delta: \PM \rw \zz$ which maps an object $(X,\a, *)$ to  $(X, \{ d_\l \ | \ 0 < \l \leq 1 \}).$
\begin{proof}
That $(X, \{ d_\l \ | \ 0 < \l \leq 1 \})$ belongs to $\zz$ goes exactly is the same way as the proof of {\bf(US), (UD), (UT)} with respect to the same continuous t-norm $*$ and {\bf(UH)}.

Next suppose $f: (X,\a) \rw (Y, \a')$ is a non-expansive map.
First observe that in case $d'_\l (f(x), f(x'))= \infty$ then for every $\gamma < \infty,$
$ \a(f(x), f(x'), \gamma) \not > 1 - \l.$
This implies that for every $\gamma < \infty, \ \a(x,x',\gamma) \not > 1 - \l.$ Hence we can conclude that $d_\l(x,x') = \infty.$

 In the finite case, applying \eqref{isoequiv} we obtain the following equivalences:
  \begin{eqnarray*}\label{morfequiv}
  f: (X,\a) \rw (Y, \a')\  \text{is non-expansive} \ &\Leftrightarrow & \forall t, \ \a(x,x',t) \leq \a'(f(x), f(x'), t )\\
&\Leftrightarrow& \forall t, \forall \l : \ \ \a(x,x', t) > 1 - \l  \Rightarrow\\
&& \a'(f(x), f(x') , t) > 1 - \l \\
&\Leftrightarrow& \forall t,  \forall \l :  \  d_\l(x,x') < t \Rightarrow  d'_\l(f(x), f(x')) < t\\
&\Leftrightarrow& \forall \l, \  d'_\l(f(x), f(x')) \leq d_\l (x,x').\\
\end{eqnarray*}
It follows that $\Delta f: (X, \{ d_\l \ | \ 0 < \l \leq 1 \}) \rw (Y, \{ d'_\l \ | \ 0 < \l \leq 1 \}) $ is levelwise non-expansive.
\end{proof}
\end{proposition}

Next we investigate the other direction.
Let $(X, \{ d_\l \ | \ 0 < \l \leq 1 \})$ be an object of $\zz$ and let $(X, \b, *)$ be defined as in \ref{reconstruction}:
$\b (x,y, \infty ) = 1$ and further for $\gamma < \infty$:
\begin{equation}\label{isomorphbeta}
\b(x,y,\gamma) = \sup \{1 - \l \ | \ d_\l(x, y) < \gamma\} = 1 - \inf  \{\l \ | \ d_\l (x, y) < \gamma\}.
\end{equation}

Again for $\gamma< \infty$ we have
\begin{equation} \label{morphi}
\b(x,y,\gamma) > 1 - \l   \Leftrightarrow  d_\l(x, y) < \gamma.
\end{equation}
As we know from the previous section, due to the properties of  $(X, \{ d_\l \ | \ 0 < \l \leq 1 \})$ listed above, $(X, \b, *)$ defines a probabilistic metric space.

\begin{proposition}\label{Phi}
This defines a concrete functor $\Phi : \zz \rw \PM$ which maps an object $(X, \{ d_\l \ | \ 0 < \l \leq 1 \})$ to $(X, \b, *)$. 
\begin{proof}
That this is indeed a functor goes as follows.\\
Let $f: (X, \{ d_\l \ | \ 0 < \l \leq 1 \}) \rw (X, \{ d'_\l \ | \ 0 < \l \leq 1 \})$ be levelwise non-expansive.
First observe that in case $\gamma = \infty,$  $\b(x,x',\infty) =  \b'(f(x), f(x'), \infty )= 1$ by definition.
 In the finite case, applying \eqref{morphi} we obtain the following equivalences:
  \begin{eqnarray*}\label{phifunctor}
  f:  (X, \{ d_\l \ | \l \}) \rw (X, \{ d'_\l \ | \l \}) \  \text{levelw. n-exp.}&\Leftrightarrow & \forall \l, \  d'_\l(f(x), f(x')) \leq d_\l (x,x')\\
&\Leftrightarrow& \forall \gamma,  \forall \l : \  d_\l(x,x') < \gamma \Rightarrow\\
&&  d'_\l(f(x), f(x')) < \gamma\\
&\Leftrightarrow& \forall \gamma, \forall \l : \ \ \b(x,x', \gamma) > 1 - \l  \Rightarrow\\
&& \b'(f(x), f(x') , \gamma) > 1 - \l )\\
&\Leftrightarrow&\b(x,x', -) \leq \b'(f(x), f(x') , -)\\
&\Leftrightarrow&\Phi(f) : (X, \b ) \rw (Y, \b' )\  \text{n-exp.} 
\end{eqnarray*}
\end{proof}
\end{proposition}

\begin{theorem}
$\Delta: \PM \rw \zz$ is an isomorphism.
\begin{proof}
We prove that $ \Phi \circ \Delta =\mathsf{id}_{\PM} $ and $ \Delta \circ \Phi = \mathsf{id}_{\zz}.$

Given a probabilistic metric space $(X, \a, *),$ applying $\Delta$ we obtain an object  $(X, \{d_\l \ | \ 0 <\l \leq 1 \}$ in $\zz$. Next applying $\Phi$ we obtain a probabilistic metric space $(X, \b, *),$ as described in \ref{delta} and \ref{Phi}. Applying \eqref{isoequiv} and \eqref{morphi}  for $\gamma < \infty$ we get 
$$
\b(x,y,\gamma) > 1 - \l  \Leftrightarrow d_\l (x,y) < \gamma \Leftrightarrow \a(x,y,\gamma) > 1-\l.
$$
Since $\b(x,y,\infty) = \a(x,y,\infty)= 1$, we have $\a = \b.$ 

For the other identity, we start with an object  $(X, \{d_\l \ | \ 0 <\l \leq 1 \})$ of $\zz$ and apply $\Phi$. Let $(X, \b, *)$ be the probabilistic metric space obtained. Applying $\Delta$ we get the object $(X, \{e_\l \ | \ 0 <\l \leq 1 \})$ from $\zz.$
Again applying \eqref{isoequiv} and \eqref{morphi}, for $\gamma < \infty$ we get
$$
e_\l(x,y) < \gamma \Leftrightarrow \b(x,y, \gamma) > 1 - \l \Leftrightarrow d_\l (x,y)  < \gamma,
$$
which implies $e_\l = d_\l$ for all $\l.$
\end{proof}

\end{theorem}

In some papers on the subject the distances involved are assumed to be finitely valued. This condition can be obtained in our setting as well.

\begin{definition}
 Define $\PM_{\text{lim}}$ the concrete full subcategory of $\PM$ consisting of those objects $(X, \a, *)$ satisfying
 \begin{equation}\label {lim1}
 \forall x,y \in X,  \lim_{\gamma \rw \infty} \a(x,y,\gamma) = 1.
 \end{equation} 
 Define $\zz_{\text{fin}}$ the full subcategory of $\zz$ consisting of those $(X, \{ d_\l \ | \ 0 < \l \leq 1 \})$ satisfying
 \begin{equation}\label{fin}
 \forall x,y \in X, \ \forall \  0 < \l \leq 1, \ d_\l (x,y) < \infty.
 \end{equation}
\end{definition}

\begin{proposition}
The functors $\Delta$ and $\Phi$ restrict to $ \Delta : \PM_{\text{lim}} \rw  \zz_{\text{fin}}$ and $\Phi :  \zz_{\text{fin}} \rw \PM_{\text{lim}}$ and hence the categories 
$\PM_{\text{lim}}$ and  $\zz_{\text{fin}}$ are isomorphic.
\begin{proof}
Let $x,y \in X$ and suppose $\lim_{\gamma \rw \infty} \a (x,y, \gamma) = 1$. This implies that for every $0 < \l \leq 1,$ the set $$ \{ 0 \leq \gamma < \infty \ | \ \a(x,y,\gamma) > 1 - \l \}$$ is non-empty and hence by \eqref{distanceassociated} $d_\l(x,y) < \infty.$

Let $x,y \in X$ and suppose that for every $0 < \l \leq 1$ the distance $d_\l(x,y) < \infty.$ Let $\l$ be arbitrary, take $\l'<\l$ and then $\gamma$ satisfying
$d_{\l'}(x,y) < \gamma.$ This implies $$ \inf\{ \rho \ | \ d_\rho(x,y) <\gamma \} \leq \l' < \l.$$
Applying \eqref{isomorphbeta} we get $\b(x,y, \gamma) > 1 - \l.$

\end{proof}
\end{proposition}

\vspace{1cm}

\vspace{1cm}
\noindent {E. Colebunders \\
Department of Mathematics and Data Science, Vrije Universiteit Brussel, Pleinlaan 2, 1050 Brussel, Belgi\"{e}\\  \emph{evacoleb@vub.be}\\
Department of Mathematics, Universiteit Antwerpen, Middelheimlaan 1, 2020 Antwerpen, Belgi\"{e}}\\
\emph{eva.colebunders@uantwerpen.be}\\
 \\
R. Lowen \\
Department of Mathematics, Universiteit Antwerpen, Middelheimlaan 1, 2020 Antwerpen, Belgi\"{e}}\\
\emph{bob.lowen@uantwerpen.be}

\end{document}